\documentclass[a4paper,12pt]{amsart}
\usepackage{amssymb,latexsym}
\def\sideremark#1{\ifvmode\leavevmode\fi\vadjust{\vbox to0pt{\vss
 \hbox to 0pt{\hskip\hsize\hskip1em
 \vbox{\hsize2cm\tiny\raggedright\pretolerance10000
 \noindent #1\hfill}\hss}\vbox to8pt{\vfil}\vss}}}%

\newtheorem{thm}{Theorem}
\newtheorem{propo}[thm]{Proposition}
\newtheorem{lem}[thm]{Lemma}

\def\ldwr{\langle d \w \rangle}

\def\cD{{\it Div}}
\def\cQ{{\cal{Q}}}
\def\cS{{\cal{S}}}
\let\z=\zeta
\def\ci{C^\infty}

\let\implies=\Rightarrow
\def\tr{{\rm tr}}

\newcommand{\comment}[1]{}

\def\aligned{\begin{array}{rl}}
\def\endaligned{\end{array}}

\newcommand{\End}{\operatorname{End}}

\def\J{{\sf J}}

\def\P{{\sf P}}
\def\R{{\sf R}}

\let\d=\delta
\let\D=\Delta
\let\e=\varepsilon

\let\l=\lambda

\let\nd=\nabla

\let\s=\sigma

\let\S=\Sigma
\let\w=\omega

\def\bbR{{\mathbb{R}}}

\def\cC{{\cal C}}

\def\cF{{\cal F}}

\def\kth{k^{\underline{{\rm th}}}}

\def\bear{\begin{array}}
\def\eear{\end{array}}

\def\beq{\begin{equation}}
\def\eeq{\end{equation}}

\def\kth{k^{\underline{{\rm th}}}}

\def\cal{\mathcal}

\renewcommand{\dfrac}[2]{{\displaystyle{\frac{#1}{#2}}}}

\newcommand{\nn}[1]{(\ref{#1})}

\renewcommand{\text}[1]{{\rm#1}}

\newcommand{\rpl}                         
{\mbox{$
\begin{picture}(12.7,8)(-.5,-1)
\put(0,0.2){$+$}
\put(4.2,2.8){\oval(8,8)[r]}
\end{picture}$}}

\newcommand{\lpl}                         
{\mbox{$
\begin{picture}(12.7,8)(-.5,-1)
\put(2,0.2){$+$}
\put(6.2,2.8){\oval(8,8)[l]}
\end{picture}$}}

\newcommand{\wh}[1]{\widehat{#1}}

\def\Scal{{\sf Scal}}
\def\Ric{{\sf Ric}}
\def\Rm{{\sf R}}
\def\Cot{{\sf C}}
\def\Wl{{\sf W}}

\newcommand{\hook}{\raisebox{-0.35ex}{\makebox[0.6em][r]
{\scriptsize $-$}}\hspace{-0.15em}\raisebox{0.25ex}{\makebox[0.4em][l]{\tiny
 $|$}}}

\def\whatami{\frac{n!}{[(n/2)!]^2}}

\newcommand{\itemi}[1]{{\rm(}#1{\rm)}}

\begin{document}
\renewcommand{\today}{}
\title{Variational status of a class of 
fully nonlinear curvature prescription problems}
\author{Thomas P.\ Branson and A.\ Rod Gover}

\address{Department of Mathematics, The 
University of Iowa, Iowa City IA
52242 USA}

\address{Department of Mathematics\\
  The University of Auckland\\
  Private Bag 92019\\
  Auckland 1\\
  New Zealand} \email{gover@math.auckland.ac.nz}

\vspace{10pt}

\renewcommand{\arraystretch}{1}
\renewcommand{\arraystretch}{1.5}

\pagestyle{myheadings}   
\dedicatory{The second named author dedicates the paper to the memory of
Thomas P. Branson (1953 - 2006)}

\begin{abstract}
Prescribing, by conformal transformation, the $k^{\rm th}$-elementary
symmetric polynomial of the Schouten tensor $\s_k(\P)$ to be constant
is a generalisation of the Yamabe problem. On compact Riemannian
$n$-manifolds we show that, for $3\leq k\leq n$, this prescription
equation is an Euler-Lagrange equation of some action if and only if
the structure is locally conformally flat.
\end{abstract}
\maketitle \markboth{Branson \& Gover}{Variational status of the 
$\s_k$-Yamabe problem}

\section{Introduction} \label{intro}

Recently there has been significant interest and progress in the study
of the so-called {\em $\s_k$-Yamabe problem} and related curvature prescription
problems, see for example \cite{V1,CGY2,GV2,LL1,GW2}. On a Riemannian
manifold $(M,g)$ (of dimension $n\geq 2$) we have the well-known
decomposition of the Riemannian curvature $\R=\Wl + \P\odot g$ where
$\Wl$ is the totally-trace-free Weyl tensor, $\P$ is the {\em Schouten
  tensor} and $\odot$ denotes the Kulkarni-Nomizu product.  For a
given fixed $k$ ($1\leq k\leq n$) the $\s_k$-Yamabe problem is to
find, within a conformal class of metrics, a metric for which
\begin{equation}\label{prescribe}
\sigma_k(\P) = \mbox{ constant },
\end{equation} 
where $\sigma_k(\P)$ is the $k^{\rm th}$-elementary symmetric
polynomial of $\P$. For $k=1$ this is the classical Yamabe problem.
Its solution by Schoen, Aubin, Trudinger, and Yamabe
(see \cite{L&P}) was a milestone in differential geometry.

We investigate here the question of whether the equation
\nn{prescribe} is variational, that is, whether it is the
Euler-Lagrange equation of some functional. This is obviously
important for the treatment of this prescription curvature problem. In
particular our study is partly motivated by \cite{nt} where it is
shown that \nn{prescribe} has a solution in settings where it is
variational. This class of settings includes the cases $k=2$ and when
$M$ is conformally flat \cite{V1,nt}.  Our main result is that on
compact Riemannian $n$-manifolds the following holds.
\begin{thm}\label{all} For $3\le k\le n$, the quantity
$\s_k(\P)$ is variational 
in a conformal class $\cC$ if and only if $\cC$ is locally flat.
\end{thm}
In the cases where $\s_{n/2}(P)$ ($n$-even) is variational, the
problem of explicitly finding an action is rather different to the
situation for other $k$. It turns out that this is related to the
problem of finding an action for the Q-curvature. This is made precise
in section section \ref{action}, and dimensions 4 and 6 are treated
explicitly.

The authors are grateful to the MSRI Berkeley, the organisers of the
Fall 2005 programme there ``Nonlinear Elliptic Equations and their
Applications'', and the organisers of the workshop ``Recent Results in
Nonlinear Elliptic Equations and their Interactions with Geometry''.
This article was conceived and largely developed during those events.
Thanks are also due to Neil Trudinger for helpful discussions.

\section{Background}\label{back}

We will say a scalar function is a {\em local scalar invariant} if it
is a natural scalar, that is, it is a quantity built polynomially from
the metric $g$ and its inverse, and the covariant derivative $\nd$ and
Riemann curvature $\Rm$ associated with $g$. (These are really the
{\em even} invariants.  On oriented manifolds one may construct
further natural invariants via the volume form. However we may ignore
this class of invariants for our current purposes.)  Suppose now (and
henceforth) that $M$ is compact.  A local scalar invariant $L$ is
{\em (conformally) variational} within a conformal class of metrics
$\cC=\{\wh{g}=e^{2\w}g\mid\w\in C^\infty(M)\}$ on a manifold $M$ if
there is a scalar valued functional (called an {\em action} or {\em
Lagrangian}) $\cF(g)$ on $\cC$ with
\begin{equation}\label{Fbul}
\cF^\bullet(g)(\w)=\int_M\w L\,dv_g\,,\qquad\mbox{all }\w\in\ci(M).
\end{equation}
Here $dv_g$ is the pseudo-Riemannian measure, and
\begin{equation}\label{bul}
\cF^\bullet(g)(\w):=\dfrac{d}{d\e}\Bigg|_{\e=0}\cF(e^{2\e\w}g).
\end{equation}
Below, we shall use this bullet notation for the conformal variation
in various contexts, sometimes
suppressing mention of the initial metric $g$.
In \nn{bul}, of course, the curve of metrics $e^{2\e\w}g$ may be replaced by
any curve with the same initial tangent $g^\bullet=2\w g$.
We stress that the property of being variational depends both on $L$ and
on the conformal class $\cC$. 

Suppose now that $g$ has Riemannian signature.
Given a variational local invariant $L$ and the corresponding functional
$\cF$, suppose that $g$ is a critical metric for 
conformal variations $\w$ with 
$\overline{\w}:=\int_M\w\,dv_g=0$.  
(These may be interpreted as volume-preserving perturbations.)  Then
the quantity in \nn{Fbul} vanishes for all $\w\in\ci(M)$ with
$\overline{\w}=0$, and this in turn implies that $L$ is constant.
Thus for variational quantities $L$, prescription of a constant value
through conformal deformation is the Euler-Lagrange equation for the
functional $\cF$.  A classic example (and the simplest example under
consideration in this paper) is the Yamabe problem of prescribing
constant scalar curvature. With $\J$ a constant multiple of the scalar
curvature, here the functional may be taken to be $\int_M\J\,dv_g\,$
if $n\ge 3$, and
$$
\frac14\int_M\log(g/g_0)(\J\,dv+(\J\,dv)_0)\qquad\mbox{when }n=2,
$$   
where $g_0$ is any choice of background
metric from the conformal class, and for example $Jdv$ means
$J^gdv_g$.  The display gives a special case of the action for the
Q-curvature; see Section \ref{action} below.

The $\kth$ elementary symmetric polynomial of an linear endomorphism
$A=(A^a{}_b)$ is 
\begin{equation}\label{defsig}
\s_k(A)=
A^{a_1}{}_{[a_1}\cdots A^{a_k}{}_{a_k]}\,,
\end{equation}
where we use abstract index notation.  In particular, an index
occurring twice, once up and once down, denotes a contraction, and
square brackets denote antisymmetrisation. Up to a nonzero constant
factor, this gives the $\l^{(n-k)}$ coefficient in $\s_n(\l\,{\rm Id}-A)$,
where $n$ is the dimension of the space on which $A$ acts, and $\s_n$
is as in \nn{defsig} for $k=n$. This follows from an easy argument in
exterior algebra, or one may use Lemma \ref{crunch}
below. Note that $\s_n(A)$ is the determinant of $A$
and, of course, $\s_k(A)$ vanishes for $k>n$.  

The Schouten tensor, mentioned above, is a trace adjustment of the
Ricci curvature (at least for dimensions $n\geq 3$)
$$
\P=\dfrac{\Ric-\J g}{n-2}\,,\qquad\mbox{where }\J:=\dfrac{\Scal}{2(n-1)}\,.
$$ Here $\Scal=\Ric^a{}_a$ is the scalar curvature. Via the metric $\P$ yields 
a section of $\End TM$  and so we may speak of
the scalar invariants $\s_k(\P)$.  For example, $\s_1(\P)=\J$, and
$\s_2(\P)=\frac12(\J^2-|\P|^2)$.  Here, to get the norm-squared of any
tensor, we contract it against against the same tensor with abstract
indices in the same order; for example $|\P|^2=\P^{ab}\P_{ab}$.

A local scalar invariant $L$ has {\em weight} $-\ell$ if uniform dilation
of the metric has the effect
$L[A^2g]=A^{-\ell}L[g]$ for all $0<A\in\bbR$.
For example, $\s_k(\P)$ has weight $-2k$.

\begin{lem}\label{vtl} Let $L$ be a weight $-\ell$ local invariant.
\newline
\itemi{i} The map
\begin{equation}\label{ado}
D:\w\mapsto\ell\w L+L^\bullet(\w)
\end{equation}
is a natural differential operator of the form $Td$ (i.e.\ with the exterior
derivative $d$ as a right composition factor).\newline
\itemi{ii} $L$ is variational in a Riemannian conformal class $\cC$ on a compact
manifold $M$ if and only if
$D$ is formally self-adjoint at all metrics in $\cC$.\newline
\itemi{iii} If $L$ is variational in a compact Riemannian $(M,\cC)$, then
\begin{equation}\label{intinvt}
\left(\int_ML\,dv_g\right)^\bullet(\w)=(n-\ell)\int_M\w L\,dv_g\qquad
{\rm in}\ \cC.
\end{equation}
In particular, if $\ell\ne n$, then 
\begin{equation}\label{easyLag}
(n-\ell)^{-1}\int_ML\,dv_g
\end{equation}
is an action for
$L$ in $\cC$.

\end{lem}

{\bf Proof}: The conformal variation the Levi-Civita connection is a
linear function of $d \omega$. This determines the conformal variation
of the scalar invariant $L$ as a differential operator on $\omega$
except for a term $-\ell \w L $ arising from the use of the inverse
metric in making contractions. So \nn{ado} is a differential operator
of the form $Td$. (For an inductive argument see, for example,
\cite{tbms}.)  This proves (i).  Now fix a conformal class $\cC$.
Since the space of metrics in $\cC$ is contractible, 
 $\int\w L\,dv_g$ is a variation of
some functional on $\cC$ if and only if the putative second variation
$$
\cS(\eta,\w)=\int\eta(L dv_g)^\bullet
=(n-\ell)\int\eta\w L dv_g+
\int\eta(D\w)dv_g
$$
is symmetric. From the extreme right end of the display it is clear
that this holds if and only if $D$ is formally self-adjoint.  This
proves (ii).
Now
$$
\left(\int_ML\,dv_g\right)^\bullet(\w)=
\cS(1,\w)=(n-\ell)\int_M\w L\,dv_g+\int(D\w)dv_g\,.
$$
But in the variational case, 
$$
\int(D\w)dv_g=\int(D1)\w\,dv_g=\int(Td1)\w\,dv_g=0,
$$
as desired for \nn{intinvt}.$\qquad\square$

In the subsequent calculations we will write $S_{|a}$, or sometimes
$\nabla S$, for the Levi-Civita covariant derivative of a tensor
$S$. Since the Levi-Civita connection is torsion-free, for a function
$\w$, $\w_a:=\w_{|a}$ is the exterior derivative of $\w$ and
$\w_{ab}=\w_{|ab}$ is symmetric.  If $u_a$ is a 1-form and $u_{a|b}$
its (Levi-Civita) covariant derivative with respect to the metric $g$
then, the covariant derivative with respect to the conformally related
metric $\hat{g}=e^{2\Upsilon} g$ ($\Upsilon \in C^\infty(M)$) is
$$
\widehat{u_{a|b}}= u_{a|b} - \Upsilon_a u_b -\Upsilon_b u_a 
+g_{ab}\Upsilon^cu_c
$$ It follows easily that the Weyl curvature $\Wl_{ab}{}^c{}_d$ is
conformally invariant, and we have the following result for the
Schouten. 
\begin{lem} \label{translem} If $\hat{g}=e^{2\Upsilon} g$ then 
$ \widehat{\P}_{ab}=\P_{ab}-\Upsilon_{ab}+\Upsilon_a\Upsilon_b
-\frac{1}{2}\Upsilon^c\Upsilon_c g_{ab} $ and so
$$
\P^\bullet(\w)=-{\rm Hess}\,\w,
$$  where 
${\mbox {\rm Hess}}\,\w$ is the covariant Hessian 
$\w_{ab}\,$.
\end{lem}

\subsection{The main constructions}\label{main}

An interesting special case of our problem concerns
$\s_3=\P^a{}_{[a}\P^b{}_b\P^c{}_{c]}\,$. (We shall often write simply
$\s_k$ to mean $\s_k(\P)$.)  We shall show:

\begin{propo}\label{three}
Let $n\ge 3$.  The quantity $\s_3(\P)$ is variational 
on a conformal class $\cC$ if and only if $\cC$ is locally flat.
\end{propo}

This is just a special case from Theorem \ref{all} above, which makes
the same statement about $\s_k$ for $3\le k\le n$.  We present the
proof of this first because it brings out most of the main issues in a
very simple setting. Here and below (via the metric) we view $\P$ as a $(1,1)$-tensor.

{\bf Proof of Proposition \ref{three}}:
Since $\P^\bullet=-{\rm Hess}\,\w$, we have
$$
(\s_3)^\bullet(\w)=
-6\w\s_3-3\w^a{}_{[a}\P^b{}_b\P^c{}_{c]}\,.
$$
By Lemma \ref{vtl}, $\s_3(\P)$ will be variational in $\cC$ if and only if
$$
D:\w\mapsto\w^a{}_{[a}\P^b{}_b\P^c{}_{c]}
$$
is formally self-adjoint at all metrics in $\cC$.
This is equivalent to the assertion that for all functions $\w,\eta$,
$$
\eta D\w-\w D\eta\in\cD,
$$
where $\cD$ is the space of exact divergences.
But 
$$
\eta D\w\in-\w^a(\P^b{}_{[b}\P^c{}_{c}\eta)_{|a]}+\cD
=-\w^a\P^b{}_{[b}\P^c{}_{c}\eta_{a]}
-\w^a(\P^b{}_{[b}\P^c{}_{c})_{|a]}\eta+\cD.
$$
The first term on the extreme right is manifestly symmetric in $\w$ and 
$\eta$; thus
$$
\eta D\w-\w D\eta\in
(\w\eta^a-\eta \w^a)(\P^b{}_{[b}\P^c{}_{c})_{|a]}+\cD.
$$
Evidently if the natural tensor $T_a:=(\P^b{}_{[b}\P^c{}_{c})_{|a]}$
vanishes identically then $D$ is formally self-adjoint.  On the other
hand if we assume that $D$ is formally self-adjoint then, for all smooth
functions $\w, \eta$, we have 
\begin{equation}\label{diff}
(\w \eta^a -\eta \w^a)T_a\in \cD
\end{equation} 
In particular we may take $\w=1$, whence $\eta^a T_a\in \cD$ which, in turn,
implies that $\eta T^a{}_{|a}\in \cD$ for any smooth function $\eta$.
Thus $T^a{}_{|a} $ vanishes identically. But this with \nn{diff} implies
that $\eta \w^aT_a\in \cD $, for all smooth functions $\w, \eta$ and
hence $T_a= (\P^b{}_{[b}\P^c{}_{c})_{|a]}$ vanishes identically.
We conclude that $\s_3(\P)$
is variational in $\cC$ if and only if $(\P^b{}_{[b}\P^c{}_{c})_{|a]}$ vanishes
identically for each metric in $\cC$.

But, performing the antisymmetrisation indicated, we find that
\begin{equation}\label{crux}
3(\P^b{}_{[b}\P^c{}_{c})_{|a]}=\P^b{}_c\Cot^c{}_{ab}\,,
\end{equation}
where $\Cot^c{}_{ab}$ is the Cotton tensor $2\P^c{}_{[a|b]}\,$. We
recall that, from the contracted Bianchi identity, this is completely
trace-free.

Now if the conformal class $\cC$ is locally flat, then $\Cot$
vanishes at each metric of the class.  Indeed, 
$$
(n-3)\Cot_{bcd}=\Wl_{abcd|}{}^a,
$$
where $\Wl$ is the Weyl tensor.
This establishes the ``if'' part of the proposition, since
$\Wl=0$ (resp.\ $\Cot=0$) is a necessary and sufficient condition
for local conformal flatness in dimension $n\ge 4$ (resp.\ $n=3$).

To establish the ``only if'' part of the proposition, 
we need to show that the vanishing of $\P^b{}_c\Cot^c{}_{ab}$
at each metric of $\cC$ implies that $\cC$ is locally flat.
Taking the divergence, we have
\begin{equation}\label{six}
0=(\P^b{}_c\Cot^c{}_{ab})_|{}^a=\P^b{}_c\Cot^c{}_{ab|}{}^a
+\P^b{}_{c|}{}^a\Cot^c{}_{ab}
=\P^b{}_c\Cot^c{}_{ab|}{}^a-\frac12|\Cot|^2
\end{equation}
at each $g\in\cC$.

A {\em conformal normal scale} at a point $p$ (the {\em node})
is a metric in our conformal
class $\cC$ in which the $\ell^{\underline{{\rm th}}}$ symmetrised covariant
derivative of $\P$ at $p$ vanishes for $\ell=0,\cdots,m$ \cite{Goadv}. 
 Here $m$ just needs
to be chosen large enough for the particular problem at hand; in our 
case $m=0$ suffices.  By \nn{six}, $|\Cot|^2$ and thus 
the tensor $\Cot$
vanish at the node in a 
conformal normal scale $g$. Thus, for every point $p\in M$, there is a normal scale $g_p\in \cC$ so that  $\Cot^{g_p}(p)=0$.  
In dimension $n=3$ $\Cot$ is conformally invariant and so this
already shows that $\cC$ is locally conformally flat. For a given $p\in M$ and normal scale $g_p$, the
freedom to vary $g$ through additional conformally related normal scales
$\hat{g}_p=e^{2\eta}g_g$ allows us to realise an arbitrary element of
$T^*_pM$ as $(d\eta)_p\,$, see Lemma \ref{translem}.  But the general conformal
change law for the Cotton tensor is
\begin{equation}\label{varcot}
\wh{g}=e^{2\w}g\ \Rightarrow\ 
\wh{\Cot}_{abc}=\Cot_{abc}+\w_d\Wl^d{}_{abc}\,.
\end{equation}
Thus $\alpha\hook\Wl$ vanishes at $p$ for all $\alpha\in T^*_pM$ and,
as a consequence, $\Wl$ vanishes at $p$.  Since $p$ was arbitrary,
$\Wl=0$; this shows $\cC$ is locally flat if $n\ge 4.\qquad\square$

For the the case $k\ge 3$, we note a technical lemma.
\begin{lem}\label{crunch} 
If $F^{a_1}{}_{b_1}\cdots{}^{a_p}{}_{b_p}$ is a tensor
field (possibly taking values in an auxiliary vector bundle)
and $p+q\le n$, then
$$
F^{a_1}{}_{[a_1}\cdots{}^{a_p}{}_{a_p}\d^{a_{p+1}}{}_{a_{p+1}}\cdots
\d^{a_{p+q}}{}_{a_{p+q}]}=cF^{a_1}{}_{[a_1}\cdots{}^{a_p}{}_{a_p]}
$$
for some nonzero constant $c$.
\end{lem}

{\bf Proof}: By induction, it suffices to prove this for $q=1$.
$\qquad\square$

{\bf Proof of theorem \ref{all}}: Proceeding in analogy with the proof
of Proposition \ref{three}, we have
$$
(\s_{k})^\bullet(\w)+2k\w\s_k
=-k\w^{a_1}{}_{[a_1}\P^{a_2}{}_{a_2}\cdots
\P^{a_k}{}_{a_k]}\,.
$$
Thus by Lemma \ref{vtl}, 
$\s_k(\P)$ is variational if and only if the operator
$$
D:\w\mapsto\w^{a_1}{}_{[a_1}\P^{a_2}{}_{a_2}\cdots
\P^{a_k}{}_{a_k]}
$$ 
is formally self-adjoint.
But if $\cD$ is the space of exact divergences,
$$
\begin{array}{rl}
(D\w)\eta&\in-\w^{a_1}(\P^{a_2}{}_{[a_2}\cdots\P^{a_k}{}_{a_k}\eta)_{|a_1]}+\cD \\
&=-\w^{a_1}\P^{a_2}{}_{[a_2}\cdots\P^{a_k}{}_{a_k}\eta_{a_1]}
-\w^{a_1}(\P^{a_2}{}_{[a_2}\cdots\P^{a_k}{}_{a_k})_{|a_1]}\eta+\cD.
\end{array}
$$

The first term on the right is manifestly symmetric in $\w$ and $\eta$.
Thus (twice) the antisymmetric part, in $\w$ and $\eta$, of $(D\w)\eta$ is
$$
(\w\eta^{a_1}-\eta\w^{a_1})
(\P^{a_2}{}_{[a_2}\cdots\P^{a_k}{}_{a_k})_{|a_1]} ~,
$$
modulo terms in $\cD$.

Thus, arguing as in Proposition \ref{three}, we see that $\s_k$ is
variational on $\cC$ if and only if the one-form
$$
T_{a_1}:=
(\P^{a_2}{}_{[a_2}\cdots\P^{a_k}{}_{a_k})_{|a_1]}=\frac{k-1}{2}\Cot^{a_2}{}_{[a_2a_1}
\P^{a_3}{}_{a_3}\cdots\P^{a_k}{}_{a_k]}
$$ vanishes at each metric of $\cC$. If $\cC$ is locally conformally
flat then, as observed above, $\Cot$ vanishes identically, and so
$T_{a}$ vanishes identically and $\s_k$ is variational.

It remains to show the implication $\Rightarrow$ of the Theorem.
Pick any metric $g\in\cC$ and write the invariant $T_a$ at any
$\wh{g}=e^{2\w}g$ using the conformal change laws \nn{varcot} and
(from Lemma \ref{translem})
$$
\wh{\P}_{ab}\in\P_{ab}-\w_{ab}+\ldwr,
$$ where $\ldwr$ is the set (ideal) of expressions containing a factor
of an (undifferentiated) $d\w$.  By an elementary scaling
argument, the $s$-homogeneous contribution to this under
$\w\mapsto\l\w$, $0<\l\in\bbR$, must vanish for each $s$.  The same is
true if we restrict to any special class of conformal factors $\w$,
provided the class concerned is invariant under this scaling.
Choose a point $p$ and let $\cS_p$ be the set of (conformal factors) $\w$
with $(d\w)_p=0$ and $({\rm Hess}\,\w)_p$ a nonzero multiple of the metric
$g$.
The $(k-3)$-homogeneous part of $T_p$ is then
$$
{\rm const}\cdot\Cot^{a_2}{}_{[a_2a_1}
\P^{a_3}{}_{a_3}\d^{a_4}{}_{a_4}\cdots\d^{a_k}{}_{a_k]}={\rm const}\cdot
\Cot^{a_2}{}_{[a_2a_1}
\P^{a_3}{}_{a_3]}={\rm const}\cdot\Cot^{a_2}{}_{a_3a_1}\P^{a_3}{}_{a_2}\,,
$$ evaluated at $p$, where ``const'' is a nonzero constant which may
vary from expression to expression, and where we have used Lemma
\ref{crunch}.  Since $p$ and $g$ were arbitrary, we have $\s_k$
variational  only if $\Cot^{c}{}_{ab}\P^{b}{}_{c}$ vanishes for
each metric in $\cC$.  This puts us in the same situation as in the
proof of Proposition \ref{three}, at the point just above equation
\nn{six}.  The rest of the proof is now identical to the argument
given there.$\qquad\square$

\subsection{The action functional for $\s_{n/2}(\P)$ and the Q-curvature}
\label{action}

From Lemma \ref{vtl} it is clear that, even when a local invariant of
weight $-n$ is known to be variational, it is a non-trivial exercise
to obtain for it a conformal primitive. The question of
getting an action in such cases is related to that of the action for
the {\em Q-curvature}, which was defined in even dimensions in
\cite{seoul}. We discuss this in the current context since it gives an
illustration of the general picture, while at the same time providing
a route to explicit action formulae for $\s_{n/2}(\P)$. We should
point out that a homotopy formula giving a primitive  for
$\s_{n/2}(\P)$ on locally conformally flat structures (and also for
$n=4$) was given in \cite{BV}.

First, note that for a local invariant of weight $-n$ to occur, within
our current framework, the dimension $n$ must be even, since all
(even) local scalar invariants have even weight.  The conformal change
law for the Q-curvature is
\begin{equation}\label{Qhat}
\wh{g}=e^{2\w}g\implies\wh{Q\,dv}=(Q+P\w)dv,
\end{equation}
where $P$ is the {\em critical GJMS operator} \cite{GJMS}, a
conformally invariant, formally self-adjoint differential operator
with principal part $\Delta^{n/2}$. Here $\Delta $ is the Laplacian
$\nabla^*\nabla$, in terms of the Levi-Civita connection $\nabla$.
Consider the two-metric action
functional
$$
\cQ(g,g_0):=\frac14\int_M\log(g/g_0)\left\{(Q\,dv)_0+Q\,dv\right\}.
$$
Now vary $g$ conformally, $g^\bullet=2\w g$, keeping $g_0$ fixed.
We have
$$
\log(g/g_0)^\bullet=2\w,\qquad(Q\,dv)_0^\bullet=0,\qquad(Q\,dv)^\bullet=(P\w)dv.
$$
As a result, 
$$
\cQ(g,g_0)^\bullet=\frac12\int\w\left\{(Q\,dv)_0+Q\,dv\right\}
+\frac14\underbrace{\int\log(g/g_0)(P\w)dv}_{=\int[P\log(g/g_0)]\w dv}.
$$
We now take advantage of \nn{Qhat} in the form
$$
[P\log(g/g_0)]dv=2\left\{Q\,dv-(Q\,dv)_0\right\}
$$
to conclude that
$$
\cQ(g,g_0)^\bullet=\int\w Q\,dv.
$$
This shows that for any $g_0\in\cC$, the functional $g\mapsto\cQ(g,g_0)$ is an action
for $Q$.
For the problem of getting an action for $\s_{n/2}(\P)$, this immediately does the 
case $n=2$, since $\s_1(\P)=\J$ is the Q-curvature in that dimension.

More generally, suppose that $n$ is even and $\cC$ is locally flat.
Since $\s_{n/2}$ is then variational (by the easy part of Theorem
\ref{all} above), Lemma \ref{vtl}(iii) implies that
$\int\s_{n/2}(\P)dv_g$ is independent of $g\in\cC$.  
By \cite{BGP}, $\s_{n/2}(\P)$
is of the form $c\cdot{\sf Pff}+\eta$, where $\eta$ is an exact
divergence, $c$ is a universal constant, and the Pfaffian term ${\sf
Pff}$ is normalised so that its integral gives the Euler
characteristic. The case of the standard sphere, in which
$\P=\J g/n=g/2$, identifies the constant $c$: by Lemma \ref{crunch},
$$
\s_{n/2}(\P)=2^{-n/2}\underbrace{\d^{a}{}_{[a}\cdots\d^{b}{}_{b]}}_{n/2\;{\rm factors}}
=2^{-n/2}\whatami,
$$
so that integrating over $S^n$, 
$$
2c=2^{-n/2}\whatami\cdot{\rm vol}(S^n).
$$ 
Since $\int Q\,dv_g$ is also independent of $g\in\cC$, \cite{BGP}
also guarantees a universal constant $q$ with $Q=q\cdot{\rm Pff}+\z$,
where $\z$ is an exact divergence.  The case of the sphere shows that
$(n-1)!{\rm vol}(S^n)=2q$, so that
$$
\s_{n/2}(\P)=\dfrac{n!}{2^{n/2}[(n/2)!]^2(n-1)!}Q+{\rm(exact\;divergence)}.
$$
Note that the constant factor on the right is $1$,
$1/4$, $1/48$ if $n=2,4,6$ respectively;
these are the cases we  work out in detail.

If $n=4$, it was shown in \cite{tbbo} that
$$
\s_2(\P)=\frac14(Q-\D\J).
$$
But $\D\J$ has action $\frac12\int\J^2dv$, so
$$
\frac14\left(\cQ(g,g_0)-\int\J^2dv\right).
$$
is an explicit action for $\s_2(\P)$ in dimension 4.
One can make this alternating in the metric pair $(g,g_0)$
(as the $\cQ$-functional is already) by using the following.
\begin{propo}\label{Sig2}
Let $n=4$.  Then 
$$
\S_2(g,g_0):=\frac14\left(\cQ(g,g_0)-\int\left\{\J^2dv-(\J^2dv)_0\right\}
\right)
$$
is an action functional for $\s_2(\P)$.
 $\S_2$ satisfies the 
{\em cocycle condition}
\begin{equation}\label{coc}
\S_2(g_2,g_0)=\S_2(g_2,g_1)+\S_2(g_1,g_0)
\end{equation}
for any $g_0,g_1,g_2\in\cC$.
\end{propo}
The last claim is a consequence of the fact that $\cQ$ satisfies a 
similar cocycle condition \cite{BrGonewQ}.

Similar ideas provide an action for $\s_3(\P)$ in dimension 6, {\em provided}
that the conformal class  $\cC$ is locally flat.
Here we have \cite{GoPet}
$$ 
\begin{array}{rl}
Q&=\D^2\J+8|\nd\P|^2+16\P^{ab}\P_{ab|c}{}^c-32\tr(\P^3)
-8\J\J_c{}^c
-16\J|\P|^2+8\J^3 \\
&=\D^2\J+8|\nd\P|^2+16\J_{ab}\P^{ab}-32\J|\P|^2+64\tr(\P^3)
+8\J\D\J+8\J^3,
\end{array}
$$
where $\tr(\P^3):=\P^a{}_b\P^b{}_c\P^c{}_a\,$.  
To pass between the two expressions, we have used the
relation
$$
\P^{ab}\P_{ab|c}{}^c=\J_{ab}\P^{ab}-\J|\P|^2+n\tr(\P^3),
$$ which holds on Cotton spaces.

Now for $m$ even and dimensions $n\neq m$ write $P_m$ for $m^{\rm
th}$-order GJMS operator of \cite{GJMS}.  We define $Q_m$ to be the
curvature quantity obtained by applying $\frac{2}{n-m}P_m$ to 1. In
particular, working modulo divergences, in dimensions $n\ge 6$ (and
with $\cC$ locally flat), from expression (23) of \cite{GoPet}
and the identity mentioned we have
$$
Q_6\in\frac12(n-6)|d\J|^2+\dfrac{(n+2)(n-2)}4\J^3-4n\J|P|^2+16\tr(\P^3)+\cD.
$$
On the other hand,
$$
6\s_3=\J^3-3\J|\P|^2+2\tr(\P^3).
$$
Thus for $n\ge 6$,
$$
48\s_3\in Q_6+48(n-6)B+\cD,
$$
where
$$
48B:=-\frac12|d\J|^2-\dfrac{n+6}4\J^3+4\J|\P|^2;
$$
in particular, in dimension 6, 
$$
48B|_{n=6}=-\frac12|d\J|^2-3\J^3+4\J|\P|^2.
$$
According to \cite{seoul}, Theorem 5.6, this establishes the following.
\begin{propo}\label{Sig3}
Let $n=6$ and $\cC$ be locally conformally flat.  Then 
$$
\S_3(g,g_0)=\frac1{48}\cQ(g,g_0)+\int\left\{B\,dv-(B\,dv)_0\right\}
$$
is an action functional for $\s_3(\P)$. 
In addition, $\S_3$ is alternating,
and satisfies a cocycle condition
analogous to \nn{coc}. \end{propo}
Once again, the cocycle condition is obvious, because $\cQ$ also has these properties.

The explanation for the functional is as follows.
In dimensions $n>6$, we have
$$
\begin{array}{l}
(\s_3dv)^\bullet\in(n-6)\w\s_3dv +\cD\cdot dv, \\
(Q_6dv)^\bullet\in(n-6)\w Q_6dv+\cD\cdot dv.
\end{array}
$$ The first is from Lemma \ref{vtl}, since $\s_3 $ is
variational (by Theorem \ref{all} as $\cC$ is locally flat).  By
construction $Q_6$ satisfies an analogue of the Yamabe equation and
this implies the second of these (details are in \cite{seoul}).  We
continue the action from \nn{easyLag}
$$
(n-6)^{-1}\int\s_3
=\frac1{48}(n-6)^{-1}\int Q+\int B
$$ to dimension 6 to get $\S_3$ as above.

Since everything is explicit, we can 
check directly (restricting everything now to dimension 6)
that this gives an action for $\s_3$; that is, that 
$$
\frac1{48}\int Q\w\,dv+\left(\int B\,dv\right)^\bullet=
\int\w\s_3\, dv.
$$
Indeed, 
$$
\begin{array}{rl}
\left(|d\J|^2\right)^\bullet+6\w|d\J|^2&\in\w\left(2\D^2\J-4\J\D\J+4|d\J|^2\right)+\cD, \\
\left(\J^3\right)^\bullet+6\w\J^3&\in 6\w\left(\J\D\J-|d\J|^2\right)+\cD \\
\left(\J|\P|^2\right)^\bullet+6\w\J|\P|^2&\in\w\left(2\J\D\J-4|d\J|^2-2\langle
{\rm Hess}\,\J,\P\rangle-2\P^{ab}\P_{ab|c}{}^c-2|\nd\P|^2\right)+\cD \\
&=\w\left(2\J\D\J-4|d\J|^2-4\J_{ab}\P^{ab}-2|\nd\P|^2
+2\J|P|^2-12\tr(\P^3)\right)+\cD.
\end{array}
$$
These formulae and $\cQ(g,g_0)^\bullet=\int\w Q$ give the Proposition.

\end{document}